\documentclass{article}

\newtheorem{theorem}{Theorem}[section]
\newtheorem{corollary}[theorem]{Corollary}

\newtheorem{lemma}[theorem]{Lemma}
\newtheorem{proposition}[theorem]{Proposition}
\newtheorem{definition}[theorem]{Definition}

\evensidemargin\oddsidemargin
\begin{document}

\author{Vadim E. Levit and Eugen Mandrescu \\
Department of Computer Science\\
Holon Academic Institute of Technology\\
52 Golomb Str., P.O. Box 305\\
Holon 58102, ISRAEL\\
\{levitv, eugen\_m\}@barley.cteh.ac.il}
\title{Bipartite graphs with uniquely restricted maximum matchings and their
corresponding greedoids}
\date{}
\maketitle

\begin{abstract}
A \textit{maximum stable set }in a graph $G$ is a stable set of maximum
size. $S$ is a \textit{local maximum stable set} of $G$, and we write $S\in
\Psi (G)$, if $S$ is a maximum stable set of the subgraph spanned by $S\cup
N(S)$, where $N(S)$ is the neighborhood of $S$. A matching $M$ is \textit{%
uniquely restricted} if its saturated vertices induce a subgraph which has a
unique perfect matching, namely $M$ itself. Nemhauser and Trotter Jr. \cite
{NemhTro}, proved that any $S\in \Psi (G)$ is a subset of a maximum stable
set of $G$. In \cite{LevMan2} we have shown that the family $\Psi (T)$ of a
forest $T$ forms a greedoid on its vertex set. In this paper we demonstrate
that for a bipartite graph $G,\Psi (G)$ is a greedoid on its vertex set if
and only if all its maximum matchings are uniquely restricted.
\end{abstract}

\section{Introduction}

Throughout this paper $G=(V,E)$ is a simple (i.e., a finite, undirected,
loopless and without multiple edges) graph with vertex set $V=V(G)$ and edge
set $E=E(G).$ If $X\subset V$, then $G[X]$ is the subgraph of $G$ spanned by 
$X$. By $G-W$ we mean the subgraph $G[V-W]$, if $W\subset V(G)$. We also
denote by $G-F$ the partial subgraph of $G$ obtained by deleting the edges
of $F$, for $F\subset E(G)$, and we write shortly $G-e$, whenever $F$ $%
=\{e\} $. If $X,Y\subset V$ are disjoint and non-empty, then by $(X,Y)$ we
mean the set $\{xy:xy\in E,x\in XA,y\in Y\}$. The \textit{neighborhood} of a
vertex $v\in V$ is the set $N(v)=\{w:w\in V$ \ \textit{and} $vw\in E\}$. If $%
\left| N(v)\right| =1$, then $v$ is a \textit{pendant vertex} of $G$; by $%
\mathrm{pend}(G)$ we designate the set of all pendant vertices of $G$. We
denote the \textit{neighborhood} of $A\subset V$ by $N_{G}(A)=\{v\in
V-A:N(v)\cap A\neq \emptyset \}$ and its \textit{closed neighborhood} by $%
N_{G}[A]=A\cup N(A)$, or shortly, $N(A)$ and $N[A]$, if no ambiguity. $%
K_{n},C_{n}$ denote respectively, the complete graph on $n\geq 1$ vertices
and the chordless cycle on $n\geq 3$ vertices. By $G=(A,B,E)$ we mean a
bipartite graph having $\{A,B\}$ as its standard bipartition.

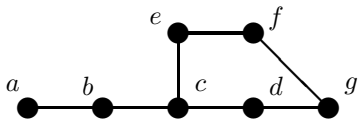
\begin{figure}[h]
\setlength{\unitlength}{1.0cm} 
\begin{picture}(5,1.5)\thicklines
  \multiput(4.5,0)(1,0){5}{\circle*{0.29}}
  \multiput(6.5,1)(1,0){2}{\circle*{0.29}}
  \put(4.5,0){\line(1,0){4}}
  \put(6.5,1){\line(1,0){1}}
  \put(6.5,0){\line(0,1){1}}
  \put(7.5,1){\line(1,-1){1}}
  \put(4.3,0.3){\makebox(0,0){$a$}}
  \put(5.3,0.3){\makebox(0,0){$b$}}
  \put(6.8,0.3){\makebox(0,0){$c$}}
  \put(7.8,0.3){\makebox(0,0){$d$}}
  \put(6.2,1.2){\makebox(0,0){$e$}}
  \put(7.8,1.2){\makebox(0,0){$f$}}
  \put(8.8,0.3){\makebox(0,0){$g$}}
\end{picture}
\caption{A graph {with diverse local maximum stable sets}.}
\label{fig10}
\end{figure}

A \textit{stable} set in $G$ is a set of pairwise non-adjacent vertices. A
stable set of maximum size will be referred to as a \textit{maximum stable
set} of $G$, and the \textit{stability number }of $G$, denoted by $\alpha
(G) $, is the cardinality of a maximum stable set in $G$. Let $\Omega (G)$
stand for the set of all maximum stable sets of $G$. A set $A\subseteq V(G)$
is a \textit{local maximum stable set} of $G$ if $A$ is a maximum stable set
in the subgraph spanned by $N[A]$, i.e., $A\in \Omega (G[N[A]])$, \cite
{LevMan2}. In the sequel, by $\Psi (G)$ we denote the set of all local
maximum stable sets of the graph $G$. For instance, any set $S\subseteq 
\mathrm{pend}(G)$ belongs to $\Psi (G)$, while the converse is not generally
true; e.g., $\{a\},\{e,d\}\in \Psi (G)$ and $\{e,d\}\cap \mathrm{pend}%
(G)=\emptyset $, where $G$ is the graph in Figure \ref{fig10}. 

Not any stable set of a graph $G$ is included in some maximum stable set of $%
G$. For example, there is no $S\in \Omega (G)$ such that $\{c,f\}\subset S$,
where $G$ is the graph depicted in Figure \ref{fig2929}. The following
theorem due to Nemhauser and Trotter Jr. \cite{NemhTro}, shows that some
special maximum stable sets can be enlarged to maximum stable sets.

\begin{theorem}
\cite{NemhTro}\label{th1} Any local maximum stable set of a graph is a
subset of a maximum stable set.
\end{theorem}

Let us notice that the converse of Theorem \ref{th1} is not generally true.
For instance, $C_{n}$, $n\geq 4$, has no proper local maximum stable set.
The graph $G$ in Figure \ref{fig10} shows another counterexample: any $S\in
\Omega (G)$ contains some local maximum stable set, but these local maximum
stable sets are of different cardinalities. As examples, $\{a,d,f\}\in
\Omega (G)$ and $\{a\},\{d,f\}\in \Psi (G)$, while for $\{b,e,g\}\in \Omega
(G)$ only $\{e,g\}\in \Psi (G)$.

In \cite{LevMan2} we have proved the following result:

\begin{theorem}
\label{th2}The family of local maximum stable sets of a forest of order at
least two forms a greedoid on its vertex set.
\end{theorem}

\begin{figure}[h]
\setlength{\unitlength}{1.0cm} 
\begin{picture}(5,1.5)\thicklines
  \multiput(4.5,0)(1,0){4}{\circle*{0.29}}
  \multiput(5.5,1)(1,0){2}{\circle*{0.29}}
  \put(4.5,0){\line(1,0){3}}
  \put(5.5,1){\line(1,0){1}}
  \multiput(5.5,0)(1,0){2}{\line(0,1){1}}
  \put(4.2,0.3){\makebox(0,0){$a$}}
  \put(5.2,0.3){\makebox(0,0){$b$}}
  \put(6.8,0.3){\makebox(0,0){$c$}}
  \put(7.8,0.3){\makebox(0,0){$d$}}
  \put(5.2,1){\makebox(0,0){$e$}}
  \put(6.8,1){\makebox(0,0){$f$}}
\end{picture}
\caption{A graph {whose }family of local maximum stable sets{\ forms a
greedoid}.}
\label{fig101}
\end{figure}
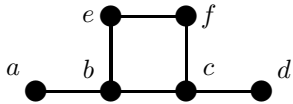

Theorem \ref{th2} is not specific for forests. For instance, the family $%
\Psi (G)$ of the graph $G$ in Figure \ref{fig101} is a greedoid.

The definition of greedoids we use in the sequel is as follows.

\begin{definition}
\cite{BjZiegler}, \cite{KorLovSch} A \textit{greedoid} is a pair $(E,%
\mathcal{F})$, where $\mathcal{F}\subseteq 2^{E}$ is a set system satisfying
the following conditions:

\setlength {\parindent}{0.0cm}(\textit{Accessibility}) for every non-empty $%
X\in \mathcal{F}$ there is an $x\in X$ such that $X-\{x\}\in \mathcal{F}$;%
\setlength
{\parindent}{3.45ex}

\setlength {\parindent}{0.0cm}(\textit{Exchange}) for $X,Y\in \mathcal{F}%
,\left| X\right| =\left| Y\right| +1$, there is an $x\in X-Y$ such that $%
Y\cup \{x\}\in \mathcal{F}$.\setlength
{\parindent}{3.45ex}
\end{definition}

Clearly, $\Omega (G)\subseteq \Psi (G)$ holds for any graph $G$. It is worth
observing that if $\Psi (G)$ is a greedoid and $S\in \Psi (G)$, $\left|
S\right| =k\geq 2$, then by accessibility property, there is a chain 
\[
\{x_{1}\}\subset \{x_{1},x_{2}\}\subset ...\subset
\{x_{1},...,x_{k-1}\}\subset \{x_{1},...,x_{k-1},x_{k}\}=S 
\]
such that $\{x_{1},x_{2},...,x_{j}\}\in \Psi (G)$, for all $j\in
\{1,...,k-1\}$. Such a chain we call an \textit{accessibility chain }for $S$%
. As an example, for $S=\{a,c,e\}\in \Psi (G)$, where $G$ is the graph in
Figure \ref{fig101}, an accessibility chain is $\{a\}\subset \{a,e\}\subset
S $.

A \textit{matching} in a graph $G=(V,E)$ is a set of edges $M\subseteq E$
having the property that no two edges of $M$ share a common vertex. We
denote the size of a \textit{maximum matching} (a matching of maximum
cardinality) by $\mu (G)$. A \textit{perfect matching} is a matching
saturating all the vertices of the graph.

Let us recall that $G$ is a \textit{K\"{o}nig-Egerv\'{a}ry graph }provided $%
\alpha (G)+\mu (G)=\left| V(G)\right| $, \cite{eger}, \cite{koen}. As a
well-known example, any bipartite graph is a K\"{o}nig-Egerv\'{a}ry graph.
Some non-bipartite K\"{o}nig-Egerv\'{a}ry graphs are presented in Figure \ref
{fig27}.

A matching $M=\{a_{i}b_{i}:a_{i},b_{i}\in V(G),1\leq i\leq k\}$ of a graph $%
G $ is called \textit{a uniquely restricted\emph{\ }matching} if $M$ is the
unique perfect matching of the subgraph $G[\{a_{i},b_{i}:1\leq i\leq k\}]$, 
\cite{GolHiLew} (first time this kind of matching appeared in \cite{HerSchn}
for bipartite graphs under the name ''\textit{constrained matching}''. Let $%
\mu _{r}(G)$ be the maximum size of a uniquely restricted matching in $G$.
Clearly, $0\leq \mu _{r}(G)\leq \mu (G)$ holds for any graph $G$. For
instance, $0=\mu _{r}(C_{2n})<n=\mu (C_{2n})$, while $\mu _{r}(C_{2n+1})=\mu
(C_{2n+1})=n$.

In this paper we characterize the bipartite graphs whose family of local
maximum stable sets are greedoids. Namely, we prove that for a bipartite
graph $G,$ the family $\Psi (G)$ is a greedoid on the vertex set of $G$ if
and only if all its maximum matchings are uniquely restricted.

Golumbic, Hirst and Lewenstein have shown in \cite{GolHiLew} that $\mu
_{r}(G)=\mu (G)$ holds when $G$ is a tree or has only odd cycles. Our
findings reveal another class of graphs enjoying this property.

\section{Preliminary results}

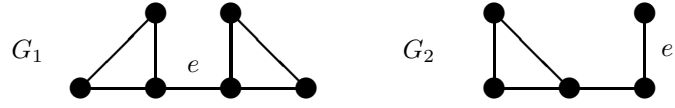
\begin{figure}[h]
\setlength{\unitlength}{1.0cm} 
\begin{picture}(5,1.5)\thicklines
  \multiput(3.5,0)(1,0){4}{\circle*{0.29}}
  \multiput(4.5,1)(1,0){2}{\circle*{0.29}}
  \put(3.5,0){\line(1,0){3}}
  \put(3.5,0){\line(1,1){1}}
  \put(4.5,0){\line(0,1){1}}
  \put(5.5,0){\line(0,1){1}}  
  \put(5.5,1){\line(1,-1){1}}
  \put(5,0.3){\makebox(0,0){$e$}}
  \put(2.8,0.5){\makebox(0,0){$G_{1}$}}

  \multiput(9,0)(1,0){3}{\circle*{0.29}}
  \multiput(9,1)(2,0){2}{\circle*{0.29}}
  \put(9,0){\line(1,0){2}}
  \put(9,0){\line(0,1){1}}
  \put(9,1){\line(1,-1){1}}
  \put(11,0){\line(0,1){1}}
  \put(11.3,0.5){\makebox(0,0){$e$}}
  \put(8,0.5){\makebox(0,0){$G_{2}$}}
\end{picture}
\caption{Non-K\"{o}nig-Egervary graphs.}
\label{fig2023}
\end{figure}

An edge $e$ of a graph $G$ is $\alpha $\textit{-critical} ($\mu $\textit{%
-critical}) if $\alpha (G)<\alpha (G-e)$ ($\mu (G)>\mu (G-e)$,
respectively). Let us observe that there is no general connection between
the $\alpha $-critical and the $\mu $-critical edges of a graph. For
instance, the edge $e$ of the graph $G_{1}$ in Figure \ref{fig2023} is $\mu $%
-critical and non-$\alpha $-critical, while the edge $e$ of the graph $G_{2}$
in the same figure is $\alpha $-critical and non-$\mu $-critical.

Nevertheless, for K\"{o}nig-Egerv\'{a}ry graphs and especially for bipartite
graphs, there is a closed relationship between these two kinds of edges.

\begin{lemma}
\cite{LevMan3}\label{critical} In a K\"{o}nig-Egerv\'{a}ry graph, $\alpha $%
-critical edges are also $\mu $-critical, and they coincide in a bipartite
graph.
\end{lemma}

In a K\"{o}nig-Egerv\'{a}ry graph, maximum matchings have a very specific
property, emphasized by the following statement:

\begin{lemma}
\cite{LevMan1}\label{match} Any maximum matching $M$ of a
K\"{o}nig-Egerv\'{a}ry graph $G$ is contained in each $(S,V(G)-\nolinebreak
S)$ and $\left| M\right| =\left| V(G)-S\right| $, where $S\in \Omega (G)$.
\end{lemma}

Clearly, not any matching of a graph is contained in a maximum matching. For
example, there is no maximum matching of the graph $G$ in Figure \ref{fig101}
that includes the matching $M=\{ab,cf\}$. Let us observe that $M$ is a
maximum matching in $G[N[\{a,f\}]],\{a,f\}$ is stable in $G$, but $%
\{a,f\}\notin \Psi (G)$. The following result shows that, under certain
conditions, a matching of a bipartite graph can be extended to a maximum
matching.

\begin{lemma}
\label{lem3}If $G$ is a bipartite graph, $\widehat{S}\in \Psi (G)$, and $%
\widehat{M}$ is a maximum matching in $G[N[\widehat{S}]]$, then there exists
a maximum matching $M$ in $G$ such that $\widehat{M}\subseteq M$.
\end{lemma}

\setlength {\parindent}{0.6cm}\textbf{Proof.} Let $W=N(\widehat{S}),$ $H=G[N[%
\widehat{S}]]$, and $S^{\prime }$ be a stable set in $G$ such that $S=%
\widehat{S}\cup S^{\prime }\in \Omega (G)$ (such $S^{\prime }$ exists
according to Theorem \ref{th1}).\setlength
{\parindent}{3.45ex} Since $H$ is bipartite and $\widehat{M}$ is a maximum
matching in $H$, it follows that 
\[
\left| \widehat{S}\right| +\left| \widehat{M}\right| =\alpha (H)+\mu
(H)=\left| V(H)\right| =\left| \widehat{S}\right| +\left| W\right| . 
\]

Let $M$ be a maximum matching in $G$. Then, by Lemma \ref{match}, $%
M\subseteq (S,V(G)-S)$, because $S\in \Omega (G)$, and 
\[
\left| M\right| =\left| V(G)-S\right| =\left| N(\bar{S})\right| +\left|
N(S^{\prime })-N(\widehat{S})\right| =\left| \widehat{M}\right| +\left|
V(G)-S-W\right| . 
\]

Let $M^{\prime }$ be the subset of $M$ containing edges having an endpoint
in $V(G)-S-W$. Since no edge joins a vertex of $\widehat{S}$ to some vertex
in $V(G)-S-W$, it follows that $M^{\prime }$ is the restriction of $M$ to $%
G[V(G)-S-W]$. Consequently, $\widehat{M}\cup M^{\prime }$ is a matching in $%
G $ that contains $\widehat{M}$, and because $\left| \widehat{M}\cup
M^{\prime }\right| =\left| \widehat{M}\right| +\left| V(G)-S-W\right|
=\left| M\right| $, we see that $\widehat{M}\cup M^{\prime }$ is a maximum
matching in $G$. \rule{2mm}{2mm}\newline

\begin{figure}[h]
\setlength{\unitlength}{1.0cm} 
\begin{picture}(5,1.5)\thicklines
  \multiput(5,0)(1,0){5}{\circle*{0.29}}
  \multiput(6,1)(2,0){2}{\circle*{0.29}}
  \put(5,0){\line(1,0){4}}
  \put(5,0){\line(1,1){1}}
  \put(6,0){\line(0,1){1}}
  \put(6,1){\line(1,-1){1}}
  \put(8,0){\line(0,1){1}}
  \put(4.8,0.3){\makebox(0,0){$a$}}
  \put(5.7,0.3){\makebox(0,0){$b$}}
  \put(5.7,1){\makebox(0,0){$c$}}
  \put(7.2,0.3){\makebox(0,0){$d$}}
  \put(8.3,1){\makebox(0,0){$e$}}
  \put(8.3,0.3){\makebox(0,0){$f$}}
  \put(9,0.3){\makebox(0,0){$g$}}
  \end{picture}
\caption{$\widehat{M}=\{ac,df\}$ is a maximum matching in $G[N[\{a,d\}]]$.}
\label{fig2929}
\end{figure}
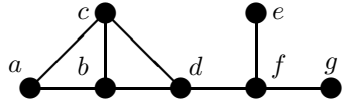

Let us notice that Lemma \ref{lem3} can not be generalized to non-bipartite
graphs. For instance, the graph $G$ presented in Figure \ref{fig2929} has $%
\widehat{S}=\{a,d\}\in \Psi (G),\widehat{M}=\{ac,df\}$ is a maximum matching
in $G[N[\widehat{S}]]$, but there is no maximum matching in $G$ that
includes $\widehat{M}$.

\begin{lemma}
\label{lem1}If $G=(A,B,E)$ is a connected bipartite graph having a unique
perfect matching, then $A\cap \mathrm{pend}(G)\neq \emptyset $ and $B\cap 
\mathrm{pend}(G)\neq \emptyset $.
\end{lemma}

\setlength {\parindent}{0.6cm}\textbf{Proof. }Let $M=\{a_{i}b_{i}:1\leq
i\leq n,a_{i}\in A,b_{i}\in B\}$ be the unique perfect matching of $G$.
Clearly, $\left| A\right| =\left| B\right| $. Suppose that $B\cap \mathrm{%
pend}(G)=\emptyset $. Hence, $\left| N(b_{i})\right| \geq 2$ for any $%
b_{i}\in B$.

Under these conditions, we shall build some cycle $C$ having half of edges
contained in $M$, and this allows us to find a new perfect matching in $G$,
which contradicts the uniqueness of $M$. We begin with the edge $a_{1}b_{1}$%
. Since $\left| N(b_{1})\right| \geq 2$, there is some $a\in
(A-\{a_{1}\})\cap N(b_{1})$, say $a_{2}$. We continue with $a_{2}b_{2}\in M$%
. Further, $N(b_{2})$ contains some $a\in (A-\{a_{2}\})$. If $a_{1}\in
N(b_{2})$, we are done, because $G[\{a_{1},a_{2},b_{1},b_{2}\}]=C_{4}$.
Otherwise, we may suppose that $a=a_{3}$, and we add to the growing cycle
the edge $a_{3}b_{3}$. Since $G$ has a finite number of vertices, after a
number of edges from $M$, we must find some edge $a_{j}b_{k}$ with $1\leq
j<k $. So, the cycle $C$ we found has 
\[
V(C)=\{a_{i},b_{i}:j\leq i\leq k\},\ E(C)=\{a_{i}b_{i}:j\leq i\leq k\}\cup
\{b_{i}a_{i+1}:j\leq i<k\}\cup \{a_{j}b_{k}\}. 
\]
Clearly, half of edges of $C$ are contained in $M$.%
\setlength
{\parindent}{3.45ex}

Similarly, we can show that also $A\cap \mathrm{pend}(G)\neq \emptyset $. 
\rule{2mm}{2mm}\newline

The following proposition presents a recursive structure of bipartite graphs
owing unique perfect matchings, which generalizes the recursive structure of
trees having perfect matching due to Fricke, Hedetniemi, Jacobs and
Trevisan, \cite{FrHedJaTre}.

\begin{proposition}
\label{prop1}$K_{2}$ is a bipartite graph, and it has a unique perfect
matching. If $G$ is a bipartite graph with a unique perfect matching, then $%
G+K_{2}$ is also a bipartite graph having a unique perfect matching.
Moreover, any bipartite graph containing a unique perfect matching can be
obtained in this way.

By $G+K_{2}$ we mean the graph comprising the disjoint union of $G$ and $%
K_{2}$, and additional edges joining at most one of endpoints of $K_{2}$ to
vertices belonging to only one color class of $G$.
\end{proposition}

\setlength {\parindent}{0.6cm}\textbf{Proof. }Let $G=(A,B,E)$ be a bipartite
graph having a unique perfect matching, say $M=\{a_{i}b_{i}:1\leq i\leq
n,a_{i}\in A,b_{i}\in B\}$. If $K_{2}=(\{x,y\},\{xy\})$, then $H=G+K_{2}$ is
also bipartite and $M\cup \{xy\}$ is a unique perfect matching in $H$, since 
$M$ was unique in $G$ and at least one of $x,y$ is pendant in $H$.%
\setlength
{\parindent}{3.45ex}

Conversely, let $G$ be a bipartite graph with a unique perfect matching. By
Lemma \ref{lem1}, it follows that $G$ has at least one pendant vertex, say $%
x $. If $y\in N(x)$, then, clearly, $G=(G-\{x,y\})+K_{2}$. \rule{2mm}{2mm}%
\newline

\section{Main results}

\begin{proposition}
\label{prop2}If $G$ is a bipartite graph of order $2n$ having a perfect
matching $M$, then $M$ is unique if and only if for some $S\in \Omega (G)$
there exists an accessibility chain.
\end{proposition}

\setlength {\parindent}{0.6cm}\textbf{Proof. }Since $\mu (G)=n$, in every
set of size greater than $n$ there exists a pair of adjacent vertices, and
hence $\alpha (G)=n$.\setlength
{\parindent}{3.45ex}

Suppose that $G$ is a bipartite graph of order $2n$ with a unique perfect
matching. We prove, by induction on $n$, that for some $S\in \Omega (G)$
there exists an accessibility chain.

For $n=2$, let $S=\{x_{1},x_{2}\}\in \Omega (G),N(S)=\{y_{1},y_{2}\}$ and $%
x_{1}y_{1},x_{2}y_{2}\in M$, where $M$ is its unique perfect matching. Then,
at least one of $x_{1},x_{2}$ is pendant, say $x_{1}$. Hence, $%
\{x_{1}\}\subset \{x_{1},x_{2}\}=S$ is an accessibility chain.

Suppose that the assertion is true for $k<n$. Let $G=(A,B,E)$ be of order $%
2n $ and $M=\{a_{i}b_{i}:1\leq i\leq n,a_{i}\in A,b_{i}\in B\}$ be its
unique perfect matching. According to Proposition \ref{prop1}, $G=H+K_{2}$.
Consequently, we may assume that: $K_{2}=(\{a_{1},b_{1}\},\{a_{1}b_{1}\})$
and $a_{1}\in \mathrm{pend}(G)$. Clearly, $H$ is a bipartite graph
containing a unique perfect matching, namely $M_{H}=M-\{a_{1}b_{1}\}$.

\emph{Case 1.} $a_{1}\in S$. Hence, $S_{n-1}=S-\{a_{1}\}\in \Omega (H)$, and
by induction hypothesis, there is a chain 
\[
\{x_{1}\}\subset \{x_{1},x_{2}\}\subset ...\subset
\{x_{1},x_{2},...,x_{n-2}\}\subset \{x_{1},x_{2},...,x_{n-1}\}=S_{n-1} 
\]
such that $\{x_{1},x_{2},...,x_{k}\}\in \Psi (H)$ for any $k\in
\{1,...,n-1\} $. Since $N(a_{1})=\{b_{1}\}$, it follows that $%
N_{G}(\{x_{1},x_{2},...,x_{k}\}\cup
\{a_{1}\})=N_{H}(\{x_{1},x_{2},...,x_{k}\})\cup \{b_{1}\}$, and therefore $%
\{x_{1},x_{2},...,x_{k}\}\cup \{a_{1}\}\in \Psi (G)$ for any $k\in
\{1,...,n-1\}$. Clearly, $\{a_{1}\}\in \Psi (G)$, and consequently, we have
the chain:

\begin{eqnarray*}
\{a_{1}\} &\subset &\{a_{1},x_{1}\}\subset \{a_{1},x_{1},x_{2}\}\subset
...\subset \{a_{1},x_{1},x_{2},...,x_{n-2}\}\subset \\
&\subset &\{a_{1},x_{1},x_{2},...,x_{n-1}\}=\{a_{1}\}\cup S_{n-1}=S,
\end{eqnarray*}
where $\{a_{1},x_{1},x_{2},...,x_{k}\}\in \Psi (G)$, for all $k\in
\{1,...,n-1\}$.

\emph{Case 2.} $b_{1}\in S$. Hence, $S_{n-1}=S-\{b_{1}\}\in \Omega (H)$ and
also $S_{n-1}\in \Psi (G)$, because $N_{G}[S_{n-1}]=A\cup B-\{a_{1},b_{1}\}$%
. By induction hypothesis, there is a chain 
\[
\{x_{1}\}\subset \{x_{1},x_{2}\}\subset ...\subset
\{x_{1},x_{2},...,x_{n-2}\}\subset \{x_{1},x_{2},...,x_{n-1}\}=S_{n-1} 
\]
such that $\{x_{1},x_{2},...,x_{k}\}\in \Psi (H)$ for any $k\in
\{1,...,n-1\} $. Since none of $a_{1},b_{1}$ is contained in $%
N_{G}(\{x_{1},x_{2},...,x_{k}\})$, it follows that $\{x_{1},x_{2},...,x_{k}%
\}\in \Psi (G)$, for any $k\in \{1,...,n-1\}$. Consequently, we have the
chain 
\[
\{x_{1}\}\subset \{x_{1},x_{2}\}\subset ...\subset
\{x_{1},x_{2},...,x_{n-1}\}=S_{n-1}\subset S_{n-1}\cup \{b_{1}\}=S, 
\]
where $\{x_{1},x_{2},...,x_{k}\}\in \Psi (G)$, for all $k\in \{1,...,n-1\}$.

Conversely, let $M=\{x_{i}y_{i}:1\leq i\leq n\}$ be a perfect matching in $G$%
, and suppose that for $S\in \Omega (G)$ there exists a chain of local
maximum stable sets 
\[
\{x_{1}\}\subset \{x_{1},x_{2}\}\subset ...\subset
\{x_{1},x_{2},...,x_{n-1}\}\subset \{x_{1},x_{2},...,x_{\alpha -1},x_{\alpha
}\}=S. 
\]

We show, by induction on $k=\left| \{x_{1},x_{2},...,x_{k}\}\right| $ that $%
H_{k}=G[N[\{x_{1},x_{2},...,x_{k}\}]]$ owns a unique perfect matching.

For $k=1$, the assertion is true, because $\{x_{1}\}\in \Psi (G)$ ensures
that $x_{1}$ is pendant, and therefore, $H_{1}=G[N[\{x_{1}\}]]$ has a unique
perfect matching, consisting of the unique edge issuing from $x_{1}$, namely 
$x_{1}y_{1}$.

Assume that $H_{k}$ has a unique perfect matching, say $M_{k}$. We may
assert that $M_{k}\subseteq M$, because $M_{k}$ is unique and included in $%
H_{k}$ and also $M$ matches $x_{1},x_{2},...,x_{k}$ onto vertices belonging
to $N(\{x_{1},x_{2},...,x_{k}\})$. Hence, $M_{k+1}=M_{k}\cup
\{x_{k+1}y_{k+1}\}$ is a maximum matching in $H_{k+1}$. If $M_{k+1}$ is not
unique in $H_{k+1}$, then there exists some $z\in
N(a_{k+1})-N[\{a_{1},a_{2},...,a_{k}\}]$ such that $z\neq y_{k+1}$.
Therefore, we infer that the set $\{x_{1},x_{2},...,x_{k}\}\cup
\{z,y_{k+1}\} $ is stable in $H_{k+1}$ and larger than $%
\{x_{1},x_{2},...,x_{k+1}\}$, which contradicts the fact that $%
\{x_{1},x_{2},...,x_{k+1}\}\in \Psi (G)$. Consequently, $M_{k+1}$ is unique
and also perfect in $H_{k+1}$. \rule{2mm}{2mm}\newline

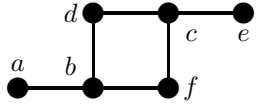
\begin{figure}[h]
\setlength{\unitlength}{1.0cm} 
\begin{picture}(5,1.5)\thicklines
  \multiput(5.5,0)(1,0){3}{\circle*{0.29}}
  \multiput(6.5,1)(1,0){3}{\circle*{0.29}}
  \put(5.5,0){\line(1,0){2}}
  \put(6.5,1){\line(1,0){2}}
  \multiput(6.5,0)(1,0){2}{\line(0,1){1}}
  \put(6.2,1){\makebox(0,0){$d$}}
  \put(7.8,0){\makebox(0,0){$f$}}
  \put(5.5,0.3){\makebox(0,0){$a$}}
  \put(6.2,0.3){\makebox(0,0){$b$}}
  \put(7.8,0.7){\makebox(0,0){$c$}}
  \put(8.5,0.7){\makebox(0,0){$e$}}
\end{picture}
\caption{{Not all }maximum matchings of a graph h{ave to be uniquely
restricted}.}
\label{fig24}
\end{figure}

If one of the maximum matchings of a bipartite graph is uniquely restricted,
this is not necessarily true for all its maximum matchings. For instance,
let us consider the bipartite graph $G$ presented in Figure \ref{fig24}. The
set of edges $M_{1}=\{ab,ce\}$ is one of uniquely restricted maximum
matchings of $G$, while $M_{2}=\{bd,cf\}$ is one of its maximum matchings,
but it is not uniquely restricted.

\begin{theorem}
\label{th10}If $G$ is a bipartite graph, then the following assertions are
equivalent:

($\mathit{i}$) there exists some $S\in \Omega (G)$ having an accessibility
chain;

($\mathit{ii}$) there exists a uniquely restricted maximum matching in $G$;

($\mathit{iii}$) each $S\in \Omega (G)$ has an accessibility chain.
\end{theorem}

\setlength {\parindent}{0.6cm}\textbf{Proof. }($\mathit{i}$) $\Rightarrow $ (%
$\mathit{ii}$) Let us consider an accessibility chain of $S\in \Omega (G)$

\[
\emptyset \subset \{x_{1}\}\subset \{x_{1},x_{2}\}\subset ...\subset
\{x_{1},x_{2},...,x_{\alpha -1}\}\subset \{x_{1},x_{2},...,x_{\alpha }\}=S, 
\]
for which we define $S_{i}=\{x_{1},x_{2},...,x_{i}\}$ and $S_{0}=\emptyset $.%
\setlength
{\parindent}{3.45ex}

Since $S_{i-1}\in \Psi (G),S_{i}=S_{i-1}\cup \{x_{i}\}\in \Psi (G)$ and $G$
is bipartite, it follows that $\left| N(x_{i})-N[S_{i-1}]\right| \leq 1$,
because otherwise, if $\{a,b\}\subset N(x_{i})-N[S_{i-1}]$, then the set $%
\{a,b\}\cup $ $S_{i-1}$ is stable in $N[S_{i-1}\cup \{x_{i}\}]$, and larger
than $S_{i}=S_{i-1}\cup \{x_{i}\}$, in contradiction with the fact that $%
S_{i}\in \Psi (G)$.

Let $\{y_{i_{j}}:1\leq j\leq \mu \}$ be such that $\{y_{i_{j}}%
\}=N(x_{i_{j}})-N[S_{i_{j}-1}]$, for all $i\in \{1,...,\alpha \}$ with $%
\left| N(x_{i})-N[S_{i-1}]\right| =1$. Hence, $M=\{x_{i_{j}}y_{i_{j}}:1\leq
j\leq \mu \}$ is a matching in $G$.

\begin{itemize}
\item  \emph{Claim 1.} $\mu =\mu (G)$, i.e., $M$ is a maximum matching in $G$%
.
\end{itemize}

Since $\left| N(x_{i})-N[S_{i-1}]\right| \leq 1$ holds for all $i\in
\{1,...,\alpha \}$, where $S_{0}=N[S_{0}]=\emptyset $, and $%
\{y_{i_{j}}\}=N(x_{i_{j}})-N[S_{i_{j}-1}]$, for all $i\in \{1,...,\alpha \}$
satisfying $\left| N(x_{i})-N[S_{i-1}]\right| =1$, it follows that $%
N(S)=\{y_{i_{j}}:1\leq j\leq \mu \}$, and this ensures that $M$ is a maximal
matching in $G$, i.e., it is impossible to add an edge to $M$ and to get a
new matching.

In addition, we have 
\[
\left| V(G)\right| =\left| N[S]\right| =\left| S\right| +\left| N(S)\right|
=\left| S\right| +\left| \{y_{i_{j}}:1\leq j\leq \mu \}\right| =\alpha
(G)+\left| M\right| , 
\]
and because $\left| V(G)\right| =\alpha \left( G\right) +\mu \left( G\right) 
$, we infer that $\left| M\right| =\mu \left( G\right) $. In other words, $M$
is a maximum matching in $G$.

\begin{itemize}
\item  \emph{Claim 2.} $M$ is a uniquely restricted maximum matching in $G$.
\end{itemize}

We use induction on $k=\left| S_{k}\right| $ to show that the restriction of 
$M$ to $H_{k}=G(N[S_{k}])$, which we denote by $M_{k}$, is a uniquely
restricted maximum matching in $H_{k}$.

For $k=1,S_{1}=\{x_{1}\}\in \Psi (G)$ and this implies that $%
N(x_{1})=\{y_{i_{1}}\}$. Clearly, $M_{1}=\{x_{1}y_{i_{1}}\}$ is a uniquely
restricted maximum matching in $H_{1}$.

Suppose that the assertion is true for all $j\leq k-1$. Let us observe that 
\[
N[S_{k}]=N[S_{k-1}]\cup (N(x_{k})-N[S_{k-1}])\cup \{x_{k}\}, 
\]
because $S_{k}=S_{k-1}\cup \{x_{k}\}$.

Further we will distinguish between two different situations depending on
the number of new vertices, which the set $N(x_{k})$ brings to the set $%
N[S_{k-1}]$.

\emph{Case 1.} $N(x_{k})-N[S_{k-1}]=\emptyset $. Hence, we obtain: 
\[
\left| V(H_{k})\right| =\left| S_{k-1}\cup \{x_{k}\}\right| +\left|
M_{k-1}\right| =\left| S_{k}\right| +\left| M_{k-1}\right| =\alpha
(H_{k})+\left| M_{k-1}\right| . 
\]
Since $\left| V(H_{k})\right| =\alpha (H_{k})+\mu \left( H_{k}\right) $, the
equality $\left| V(H_{k})\right| =\alpha (H_{k})+\left| M_{k-1}\right| $
ensures that $M_{k-1}$ is a maximum matching of $H_{k}$. Therefore, $M_{k-1}$
is a uniquely restricted maximum matching in $H_{k}$.

\emph{Case 2.} $N(x_{k})-N[S_{k-1}]=\{y_{i_{k}}\}$. Then we have: 
\[
\left| V(H_{k})\right| =\left| S_{k-1}\cup \{x_{k}\}\right| +\left|
M_{k-1}\cup \{x_{k}y_{i_{k}}\}\right| =\left| S_{k}\right| +\left|
M_{k}\right| =\alpha (H_{k})+\left| M_{k}\right| , 
\]
and this assures that $M_{k}=M_{k-1}\cup \{x_{k}y_{i_{k}}\}$ is a maximum
matching in $H_{k}$. The edge $e=x_{k}y_{i_{k}}$ is $\alpha $-critical in $%
H_{k}$, because $\{y_{i_{k}}\}=N(x_{k})-N[S_{k-1}]$. According to Lemma \ref
{critical}, $e$ is also $\mu $-critical in $H_{k}$. Therefore, any maximum
matching of $H_{k}$ contains $e$, and since $M_{k}=M_{k-1}\cup \{e\}$ and $%
M_{k-1}$ is a uniquely restricted maximum matching in $H_{k-1}=H_{k}-%
\{x_{k},y_{i_{k}}\}$, it follows that $M_{k}$ is a uniquely restricted
maximum matching in $H_{k}$.

($\mathit{ii}$) $\Rightarrow $ ($\mathit{iii}$) Let $M$ be a uniquely
restricted maximum matching in $G$. According to Lemma \ref{match}, $%
M\subseteq (S,V(G)-S)$ and $\left| M\right| =\left| V(G)-S\right| =\mu (G)$.
Therefore, $M$ is a unique perfect matching in $H=G[N[S_{\mu }]]$, where 
\[
S_{\mu }=\{x:x\in S,x\ is\ an\ endpoint\ of\ an\ edge\ in\ M\}. 
\]
It is clear that $S_{\mu }$ is a maximum stable set in $H$, because $%
N(S_{\mu })=V(G)-S$ and $S_{\mu }$ is stable. In other words, $S_{\mu }\in
\Psi (G)$. Since $H$ is bipartite and $M$ is its unique perfect matching,
Proposition \ref{prop2} implies that there exists a chain 
\[
\{x_{1}\}\subset \{x_{1},x_{2}\}\subset ...\subset \{x_{1},x_{2},...,x_{\mu
-1}\}\subset \{x_{1},x_{2},...,x_{\mu -1},x_{\mu }\}=S_{\mu }, 
\]
such that all $S_{k}=\{x_{1},x_{2},...,x_{k}\},1\leq k\leq \mu $ are local
maximum stable sets in $H$. The equality $N_{H}[S_{k}]=N_{G}[S_{k}]$
explains why $S_{k}\in \Psi (G)$ for all $k\in \{1,...,\mu (G)\}$. Let now $%
x\in S-S_{\mu }$. Then $N(x)\subseteq V(G)-S$, and therefore, $N(S_{\mu
}\cup \{x\})=V(G)-S$. Since $S_{\mu }$ is a maximum stable set in $H$ and $%
S_{\mu }\cup \{x\}$ is stable in $H\cup \{x\}=G[N[S_{\mu }\cup \{x\}]]$, we
get that $S_{\mu }\cup \{x\}$ is a maximum stable set in $H\cup \{x\}$,
i.e., $S_{\mu +1}=S_{\mu }\cup \{x\}\in \Psi (G)$. If there still exists
some $y\in S-S_{\mu +1}$, in the same manner as above we infer that $S_{\mu
+2}=S_{\mu +1}\cup \{y\}\in \Psi (G)$.

In such a way we build the following accessibility chain 
\[
\{x_{1}\}\subset \{x_{1},x_{2}\}\subset ...\subset \{x_{1},x_{2},...,x_{\mu
}\}\subset S_{\mu +1}\subset S_{\mu +1}\subset ...\subset S_{\alpha }=S. 
\]

Clearly,\textbf{\ }($\mathit{iii}$) $\Rightarrow $ ($\mathit{i}$), and this
completes the proof. \rule{2mm}{2mm}\newline

\begin{figure}[h]
\setlength{\unitlength}{1.0cm} 
\begin{picture}(5,1.5)\thicklines
  \multiput(4.5,0)(1,0){4}{\circle*{0.29}}
  \multiput(5.5,1)(1,0){4}{\circle*{0.29}}
  \put(4.5,0){\line(1,0){3}}
  \put(6.5,1){\line(1,0){2}}
  \multiput(5.5,0)(1,0){3}{\line(0,1){1}}
  \put(6.2,1){\makebox(0,0){$d$}}
  \put(7.8,0){\makebox(0,0){$f$}}
  \put(4.5,0.3){\makebox(0,0){$a$}}
  \put(5.2,0.3){\makebox(0,0){$b$}}
  \put(6.2,0.3){\makebox(0,0){$e$}}
  \put(5.2,1){\makebox(0,0){$c$}}
  \put(7.8,0.7){\makebox(0,0){$g$}}
  \put(8.5,0.7){\makebox(0,0){$h$}}
\end{picture}
\caption{{The chain of uniquely restricted }matchings is :$%
\{hg\},\{hg,de\},\{hg,de,cb\}$ .}
\label{fig2444}
\end{figure}
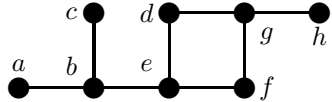

As an example of the process of building a uniquely restricted maximum
matching with the help of an accessibility chain, let us consider the
bipartite graph $G$ presented in Figure \ref{fig2444}. The accessibility
chain 
\[
\{h\}\subset \{h,d\}\subset \{h,d,f\}\subset \{h,d,f,c\}\subset
\{h,d,f,c,a\}\in \Psi (G) 
\]
gives rise to the uniquely restricted maximum matching $M=\{hg,de,cb\}$.
Notice that $\Psi (G)$ is not a greedoid, because $\{d,f\}\in \Psi (G)$,
while $\{d\},\{f\}\notin \Psi (G)$. 

The following theorem will show us another reason, why the family $\Psi (G)$
of the graph $G$ from Figure \ref{fig2444} is not a greedoid, namely $%
\{bc,de,fg\}$ is a maximum matching, but not uniquely restricted.

\begin{theorem}
\label{th11}If $G$ is a bipartite graph, then $\Psi (G)$ is a greedoid if
and only if all its maximum matchings are uniquely restricted.
\end{theorem}

\setlength {\parindent}{0.6cm}\textbf{Proof.} Assume that $\Psi (G)$ is a
greedoid. Let $M$ be a maximum matching in $G$. According to Lemma \ref
{match}, we have that $M\subseteq (S,V(G)-S)$ and $\left| M\right| =\left|
V(G)-S\right| $ for any $S\in \Omega (G)$. Let $S_{\mu }$ contain the
vertices of some $S\in \Omega (G)$ matched by $M$ with the vertices of $%
V(G)-S$. Since $M$ is a perfect matching in $G[N[S_{\mu }]]$ and $\left|
S_{\mu }\right| =\left| M\right| $, it follows that $S_{\mu }$ is a maximum
stable set in $G[N[S_{\mu }]]$, i.e., $S_{\mu }\in \Psi (G)$. Hence, there
exists an accessibility chain of the following structure: 
\[
\{x_{1}\}\subset \{x_{1},x_{2}\}\subset ...\subset \{x_{1},x_{2},...,x_{\mu
}\}=S_{\mu }\subset S_{\mu }\cup \{x_{\mu +1}\}\subset ...\subset S.
\]
While the existence of the first part of this chain, i.e., $%
\{x_{1}\},\{x_{1},x_{2}\},...,\{x_{1},x_{2},...,x_{\mu }\}$, is based on the
accessibility property of the family $\Psi (G)$, the existence of the second
part of the same chain, namely $S_{\mu },S_{\mu }\cup \{x_{\mu +1}\},...,S$,
stems from the exchange property of $\Psi (G)$. Now, according to
Proposition \ref{prop2}, we may conclude that the perfect matching $M$ is
unique in $G[N[S_{\mu }]]$. Hence, $M$ is a uniquely restricted maximum
matching in $G$.\setlength
{\parindent}{3.45ex}

Conversely, suppose that all maximum matchings of $G$ are uniquely
restricted. Let $\widehat{S}\in \Psi (G),H=G[N[\widehat{S}]]$, and $\widehat{%
M}$ be a maximum matching in $H$. The graph $H$ is bipartite as a subgraph
of a bipartite graph. By Lemma \ref{lem3}, there exists a maximum matching
in $G$, say $M$, such that $\widehat{M}\subseteq M$. Since $M$ is uniquely
restricted in $G$, it follows that $\widehat{M}$ is uniquely restricted in $%
H $. According to Theorem \ref{th10}, there exists an accessibility chain of 
$\widehat{S}$ in $H$%
\[
S_{1}\subset S_{2}\subset ...\subset S_{q-1}\subset S_{q}=\widehat{S}. 
\]
Since $N_{H}[S_{k}]=N_{G}[S_{k}]$, we infer that $S_{k}\in \Psi (G)$, for
any $k\in \{1,...,q\}$.

To complete the proof, we have to show\emph{\ }that, in addition to the
accessibility property, $\Psi (G)$ satisfies also the exchange property.

Let $X,Y\in $ $\Psi (G)$ and $\left| Y\right| =\left| X\right| +1=m+1$.
Hence, there is an accessibility chain 
\[
\{y_{1}\}\subset \{y_{1},y_{2}\}\subset ...\subset
\{y_{1},...,y_{m}\}\subset \{y_{1},...,y_{m},y_{m+1}\}=Y. 
\]
Since $Y$ is stable, $X\in $ $\Psi (G)$, and $\left| X\right| <\left|
Y\right| $, it follows that there exists some $y\in Y-X$, such that $y\notin
N[X]$. Let $M_{X}$ be a maximum matching in $H=G[N[X]]$. Since $H$ is
bipartite, $X$ is a maximum stable set in $H$, and $M_{X}$ is a maximum
matching in $H$, it follows that 
\[
\left| X\right| +\left| M_{X}\right| =\left| N[X]\right| =\left| X\right|
+\left| N(X)\right| ,\ i.e.,\left| M_{X}\right| =\left| N(X)\right| . 
\]
Let $y_{k+1}\in Y$ be the first vertex in $Y$ satisfying the conditions: $%
y_{1},...,y_{k}\in N[X]$ and $y_{k+1}\notin N[X]$. Since $%
\{y_{1},...,y_{k}\} $ is stable in $N[X]$, there is $\{x_{1},...,x_{k}\}%
\subseteq X$ such that for any $i\in \{1,...,k\}$ either $x_{i}=y_{i}$ or $%
x_{i}y_{i}\in M_{X}$.

Now we show that $X\cup \{y_{k+1}\}\in \Psi (G)$.

\emph{Case 1}. $N[X\cup \{y_{k+1}\}]=N[X]\cup \{y_{k+1}\}$. Clearly, $X\cup
\{y_{k+1}\}$ is stable in $G(N[X\cup \{y_{k+1}\}])$ and $\left| X\cup
\{y_{k+1}\}\right| =\left| X\right| +1$ ensures that $X\cup \{y_{k+1}\}\in
\Psi (G)$, because $X\in $ $\Psi (G)$ too.

\emph{Case 2}. $N[X\cup \{y_{k+1}\}]\neq N[X]\cup \{y_{k+1}\}$. Suppose
there are $a,b\in N(y_{k+1})-N[X]$. Hence, it follows that $%
\{a,b,x_{1},...,x_{k}\}$ is a stable set included in $N[\{y_{1},...,y_{k+1}%
\}]$ and larger than $\{y_{1},...,y_{k+1}\}$, in contradiction with the fact
that $\{y_{1},...,y_{k+1}\}\in \Psi (G)$. Therefore, there exists a unique $%
a\in N(y_{k+1})-N[X]$. Consequently, 
\[
N[X\cup \{y_{k+1}\}]=N[X]\cup N[y_{k+1}]=N[X]\cup \{a,y_{k+1}\}, 
\]
and since $ay_{k+1}\in E(G)$, we obtain that $X\cup \{y_{k+1}\}$ is a
maximum stable set in $G[N[X\cup \{y_{k+1}\}]]$, i.e., $X\cup \{y_{k+1}\}\in
\Psi (G)$. \rule{2mm}{2mm}\newline

As an immediate consequence of Theorem \ref{th11}, we obtain the following:

\begin{corollary}
\label{cor1}For any bipartite graph $G$ having a perfect matching, $\Psi (G)$
is a greedoid if and only if $G$ has a unique perfect matching.
\end{corollary}

\begin{figure}[h]
\setlength{\unitlength}{1.0cm} 
\begin{picture}(5,1.3)\thicklines
  \multiput(2.5,0)(1,0){4}{\circle*{0.29}}
  \multiput(3.5,1)(1,0){2}{\circle*{0.29}}
  \put(2.5,0){\line(1,0){3}}
  \put(3.5,0){\line(0,1){1}}
  \put(3.5,1){\line(1,0){1}}
  \put(4.5,1){\line(1,-1){1}}
  \put(4.5,1.3){\makebox(0,0){$u$}}
  \put(4.4,0.3){\makebox(0,0){$v$}}
  \put(3,0.3){\makebox(0,0){$e$}}
  \put(1.5,0.5){\makebox(0,0){$C_{5}+e$}}

  \put(7.5,0.5){\makebox(0,0){$C_{5}+3e$}}
  \multiput(8.5,0)(1,0){4}{\circle*{0.29}}
  \multiput(8.5,1)(1,0){4}{\circle*{0.29}}
  \put(8.5,0){\line(1,0){3}}
  \put(8.5,0){\line(1,1){1}}
  \put(8.5,1){\line(1,0){3}}
  \put(10.5,0){\line(0,1){1}}
  \put(9,1.3){\makebox(0,0){$e$}}
  \put(11,1.3){\makebox(0,0){$e$}}
  \put(11,0.3){\makebox(0,0){$e$}}
 \end{picture}
\caption{Non-bipartite graphs {with unique perfect matchings}.}
\label{fig27}
\end{figure}
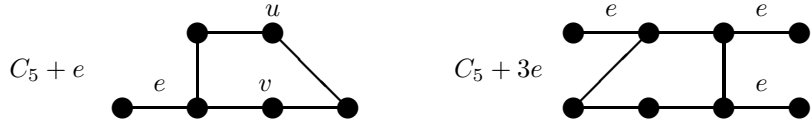

Corollary \ref{cor1} and, consequently, Theorem \ref{th11} are not valid for
non-bipartite graphs. For example, the graph $C_{5}+e$ in Figure \ref{fig27}
is a non-bipartite graph having only uniquely restricted maximum matchings,
(in fact, it has a unique perfect matching), but $\Psi (C_{5}+e)$ is not a
greedoid, because $\{u,v\}\in \Psi (C_{5}+e)$, while $\{u\},\{v\}\notin \Psi
(C_{5}+e)$.

However, there are non-bipartite graphs with a unique perfect matching,
whose $\Psi (G)$ is a greedoid. For instance, while the graph $C_{5}+3e$ in
Figure \ref{fig27} is a non-bipartite graph with a unique perfect matching,
the family $\Psi (C_{5}+3e)$ is a greedoid.

\begin{figure}[h]
\setlength{\unitlength}{1.0cm} 
\begin{picture}(5,1.5)\thicklines
  \multiput(3,0)(1,0){4}{\circle*{0.29}}
  \multiput(3,1)(1,0){3}{\circle*{0.29}}
  \put(4,0){\line(-1,1){1}}
  \put(3,0){\line(1,0){3}}
  \put(4,1){\line(1,0){1}}
  \multiput(4,0)(1,0){2}{\line(0,1){1}}
  \put(2.5,0.5){\makebox(0,0){$G_{1}$}}

  \multiput(8,0)(1,0){4}{\circle*{0.29}}
  \multiput(9,1)(1,0){3}{\circle*{0.29}}
  \put(8,0){\line(1,0){3}}
  \put(9,1){\line(1,0){2}}
  \put(9,0){\line(1,1){1}}
  \put(10,0){\line(0,1){1}}
  \put(7.5,0.5){\makebox(0,0){$G_{2}$}}
\end{picture}
\caption{$\Psi (G_{1})$ and $\Psi (G_{2})$ {form greedoids, but only }${G}_{{%
1}}${\ }is a bipartite graph.}
\label{fig20}
\end{figure}
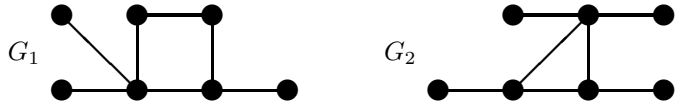

Let us also notice that there exist both bipartite and non-bipartite graphs
without a perfect matching whose family of local maximum stable sets is a
greedoid. For instance, neither $G_{1}$ nor $G_{2}$ in Figure \ref{fig20}
has a perfect matching, $G_{1}$ is bipartite, and $\Psi (G_{1}),\Psi (G_{2})$
are greedoids. 

Since any forest, by definition, has no cycles, the following Lemma \ref
{lem4} ensures that all matchings of a forest are uniquely restricted.

\begin{lemma}
\label{lem4}\cite{LevMan4} If a bipartite graph has two perfect matchings $%
M_{1}$ and $M_{2}$, then any of its vertices, from which are issuing edges
contained in $M_{1}$ and $M_{2}$, respectively, belongs to some cycle that
is alternating with respect to at least one of $M_{1}$, $M_{2}$.
\end{lemma}

It is also interesting to note that Golumbic, Hirst and Lewenstein have
proved the following generalization of Lemma \ref{lem4}.

\begin{theorem}
\label{th4}\cite{GolHiLew} A matching $M$ in a graph $G$ is uniquely
restricted if and only if there is no even-length cycle with edges
alternating between matched and non-matched edges.
\end{theorem}

Now restricting Theorem \ref{th11} to forests we immediately obtain that the
family of local maximum stable sets of a forest forms a greedoid on its
vertex set, which gives a new proof of the main finding from \cite{LevMan2},
namely Theorem \ref{th2}.

\section{Conclusions}

We have shown that to have all maximum matchings uniquely restricted is
necessary and sufficient for a bipartite graph $G$ to enjoy the property
that $\Psi (G)$ is a greedoid. We have also described all the bipartite
graphs having a unique perfect matching, or in other words, all bipartite
graphs having a perfect matching and whose $\Psi (G)$ is a greedoid. It
seems to be interesting to describe a recursive structure of all bipartite
graphs whose $\Psi (G)$ is a greedoid.

A linear time algorithm to decide whether a matching in a bipartite graph is
uniquely restricted is presented in \cite{GolHiLew}. It is also shown there
that the problem of finding a maximum uniquely restricted matching is 
\textbf{NP}-complete for bipartite graphs. These results motivate us to
propose another open problem, namely: how to recognize bipartite graphs
whose $\Psi (G)$ is a greedoid?

\end{document}